
\input amstex
\documentstyle{amsppt}
\magnification=\magstep1
\NoRunningHeads
\topmatter
\title Minimal idempotents\\
of twisted group algebras of cyclic 2-groups
\endtitle
\author Todor Zh. Mollov and Nako A. Nachev
\endauthor
\address Department of Algebra, Plovdiv University,
4000 Plovdiv, Bulgaria \endaddress
\email mollov\@ uni-plovdiv.bg and
nachev\@ uni-plovdiv.bg \endemail
\thanks Supported by the Fund ``NIMP''
of Plovdiv University under Contract 14/97. \endthanks

\abstract
For a field $K$ of characteristic different from 2, we
find the explicit form of the minimal idempotents of
the twisted group algebra $K_t\langle g\rangle$
of a cyclic 2-group $\langle g\rangle$ over $K$.
\endabstract
\endtopmatter
\document
\baselineskip 16pt

\head 1. Introduction
\endhead

Many investigations of a twisted group algebra $K_tG$ are based
on the explicit form of its minimal idempotents.
When $\langle g\rangle$ is a cyclic $p$-group, $p$ is an odd prime and
$K$ is a field of characteristic different from $p$, Nachev and Mollov
\cite{4} have found an exact formula for the minimal idempotent
of $K_t\langle g\rangle$. In \cite{1} the explicit form of the
minimal idempotents of $K_t\langle g\rangle$ has been found when
$\langle g\rangle$ is a cyclic 2-group and $K$ is a field of the
second kind with respect to 2 (and of characteristic different from 2).
In the present paper we generalize the results from \cite{1}
and find an exact formula for the minimal idempotents of
$K_t\langle g\rangle$ for $p=2$, i.e. when $\langle g\rangle$ is
a cyclic 2-group and $K$ is an arbitrary field of characteristic
different from 2.

The paper is organized as follows. In Section 1 we give some
notation, definitions and preliminary results and in Section 2
-- some preliminary results for fields $K$ of characteristic
different from 2. In Section 3 we find an explicit form of the minimal
idempotents of the twisted group algebra $K_t\langle g\rangle$ of
an arbitrary cyclic 2-group $\langle g\rangle$ over the same field $K$.

\head 1. Notation, definitions and preliminary results
\endhead

Let $K$ be a field of characteristic different
from a prime $p$ and let $\bar K$ be the algebraic closure of $K$. We
denote by $\varepsilon_n$ a $p^n$-th primitive root of 1 in
$\bar K$. The field $K$ is called a field of the first kind with
respect to $p$ [5, p. 684], if
$K(\varepsilon_i)\not= K(\varepsilon_2)$ for some $i > 2$.
Otherwise $K$ is called a field of the second kind with respect to $p$.
Let $K^{\ast}$ be the multiplicative group of $K$ and let $G$ be a
multiplicative group. A twisted group algebra $K_tG$ of $G$ over
$K$  [5, p. 13] is an associative $K$-algebra with basis
$\{\bar x\mid x\in G\}$ and multiplication defined on the basis by
$$
\bar x\bar y = \gamma(x,y)\overline{xy},\, \gamma(x,y)\in K^{\ast}.
$$
If the cyclic group $\langle g\rangle$ is of order $n$
then $\bar g^{^n} = a$ for some $a\in K^{\ast}$.
Obviously this equality determines the twisted
group algebra $K_t\langle g\rangle$ and we shall say that
the equality defines $K_t\langle g\rangle$. We define
$$
K^n =\{a^n\mid a\in K\},\,n\in{\Bbb N}.
$$

In this paper we assume that $p\not= 2$, i.e. $K$ is a field of
characteristic different from 2 and $\varepsilon_n$ is a primitive
$2^n$-th root of 1 in $\bar K$.

Let $K_2$ be the 2-component of the torsion subgroup $tK^{\ast}$
of $K^{\ast}$, i.e. $K_2$ is the Sylow 2-subgroup of $K^{\ast}$.

In order to see that some special binomials are indecomposable over
the field $K$ we shall use the following theorem.

\proclaim{Theorem A \cite{2, {\rm Theorem 16.6, p. 225}}}
The binomial $x^n - a$, $a\in K$, is irreducible over $K$ if and only
if $a\not\in K^p$ for all primes $p$ dividing $n$ and
$a\not\in -4K^4$ whenever $4\vert n$.
\endproclaim

As an immediate consequence of Theorem A one obtains the following
lemma.

\proclaim{Lemma B \cite{1, {\rm Lemma 2.7}}} The polynomial
$f(x) = x^{2^n} - \alpha$, $\alpha\in K(\varepsilon_2)$, $n\in{\Bbb N}$,
is irreducible over $K(\varepsilon_2)$ if and only if $\alpha\not\in
K(\varepsilon_2)^2$.
\endproclaim

\head 2. Some results for fields of characteristic different from 2
\endhead

\proclaim{Lemma 1} Exactly one of the following cases holds
for the field $K$:

{\rm (A)} $K=K(\varepsilon_2)$ and $\varepsilon_m\in K$
for every natural $m$;

{\rm (B)} $K=K(\varepsilon_2)$ and there exists a unique
natural $m$ such that $\varepsilon_m\in K$,
$\varepsilon_{m+1}\not\in K$ and $m\geq 2$;

{\rm (C)} $K\not= K(\varepsilon_2)$ and
$\varepsilon_m+\varepsilon_m^{-1}\in K$
for every natural $m$;

{\rm (D)} $K\not= K(\varepsilon_2)$ and there exists a unique
natural $m\geq 2$ such that
$\varepsilon_m+\varepsilon_m^{-1}\in K$,
$\varepsilon_{m+1}+\varepsilon_{m+1}^{-1}\not\in K$ and
$\varepsilon_{m+1}-\varepsilon_{m+1}^{-1}\not\in K$;

{\rm (E)} $K\not= K(\varepsilon_2)$ and there exists a unique
natural $m\geq 3$ such that
$\varepsilon_m-\varepsilon_m^{-1}\in K$.
\endproclaim

\demo{Proof} Since $K$ is a field, its Sylow 2-subgroup $K_2$ is
either cyclic or a group of type $2^{\infty}$. If $K_2$ is of
type $2^{\infty}$, then we are in the case (A). If $K_2$ is cyclic
of order $2^m$ and $K=K(\varepsilon_2)$, then we have the case (B).

Now, let $K\not=K(\varepsilon_2)$. If $\varepsilon_m\in K(\varepsilon_2)$
and $m\geq 2$, then the minimal polynomial of $\varepsilon_m$
over $K$ is of second degree. Then $\varepsilon_m$ has a unique
conjugate over $K$ element $\bar\varepsilon_m$
different from $\varepsilon_m$ and
$\varepsilon_m\bar\varepsilon_m\in\{1,-1\}$. Hence
$\bar\varepsilon_m=\varepsilon_m^{-1}$ or
$\bar\varepsilon_m=-\varepsilon_m^{-1}$.

Let the Sylow 2-subgroup $K(\varepsilon_2)_{_2}$ of
$K(\varepsilon_2)^{\ast}$ be of type $2^{\infty}$. Then for every
natural $m$ we have $\varepsilon_m\in K(\varepsilon_2)$ and
$\varepsilon_{m+1}\in K(\varepsilon_2)$. From here we obtain
$$
\varepsilon_m+\varepsilon_m^{-1}=
(\varepsilon_{m+1}\pm\varepsilon_{m+1}^{-1})^2\mp 2=
(\varepsilon_{m+1}+\bar\varepsilon_{m+1})^2\mp 2\in K.
$$
Hence we are in the case (C) and $K$ is a field of the second kind
with respect to 2.

Now, let $K(\varepsilon_2)_{_2}$ be cyclic of order $2^m$.
Clearly $m\geq 2$. If we assume that
$\varepsilon_{m+1}+\varepsilon_{m+1}^{-1}\in K$ or
$\varepsilon_{m+1}-\varepsilon_{m+1}^{-1}\in K$, we obtain that
$\varepsilon_{m+1}\in K(\varepsilon_2)_{_2}$ which is in contradiction
with the order of $K(\varepsilon_2)_{_2}$. Hence
$\varepsilon_{m+1}+\varepsilon_{m+1}^{-1}\not\in K$ and
$\varepsilon_{m+1}-\varepsilon_{m+1}^{-1}\not\in K$. Further,
as we have already mentioned, there are two possibilities for
$\bar\varepsilon_m$: either
$\bar\varepsilon_m=\varepsilon_m^{-1}$ or
$\bar\varepsilon_m=-\varepsilon_m^{-1}$. In the first case we
are in (D) and in the second in (E).

The uniqueness of the integer $m$ in the cases (B), (D) and (E)
follows from the fact, that it is determined uniquely from the order
of $K(\varepsilon_2)_{_2}$. In the case (E) we have $m\geq 3$
because for $m=2$ the condition
$\varepsilon_m-\varepsilon_m^{-1}\in K$ does not hold.
The lemma is established.
\enddemo

We want to mention that the condition
$\varepsilon_m-\varepsilon_m^{-1}\in K$ in the case (E) of the lemma
is sufficiently strong, since implies
$\varepsilon_{m+1}+\varepsilon_{m+1}^{-1}\not\in K$ and
$\varepsilon_{m+1}-\varepsilon_{m+1}^{-1}\not\in K$.
If this condition is satisfied for some $m\geq 3$, then
it does not hold for any other $m\geq 3$. The situation with
the condition
$\varepsilon_m-\varepsilon_m^{-1}\in K$ in the case (D)
is different. It is not sufficiently strong because does not imply
$\varepsilon_{m+1}+\varepsilon_{m+1}^{-1}\not\in K$ and
$\varepsilon_{m+1}-\varepsilon_{m+1}^{-1}\not\in K$.
We shall also mention that if $K$ is of the first kind with
respect to 2 and $K\not= K(\varepsilon_2)$ (i.e. in the cases
(D) and (E) of Lemma 1), then for every $l\in [2,m)\cup\{2\}$
we have $\varepsilon_l+\varepsilon_l^{-1}\in K$.
In the same two cases and in the case (B) for $l>m$ we have
$\varepsilon_{l+1}+\varepsilon_{l+1}^{-1}\not\in K$ and
$\varepsilon_{l+1}-\varepsilon_{l+1}^{-1}\not\in K$.

One sees from Lemma 1 that the class of all field of characteristic
different from 2 can be separated in 5 types which we denote
respectively with A, B, C, D, E as listed in the statement of the lemma.

In the cases (A) and (C) the field $K$ is of the second kind
with respect to 2 and in the other cases -- of the first kind.
It is not difficult to see that in the cases (B), (D) and (E), i.e.
if $K$ is of the first kind with respect to 2, then the integer $m$,
determined in Lemma 1, coincides with the constant of the field $K$
with respect to 2 (for the definition of this constant see \cite{3}).
Besides, we want to mention that this constant is equal to
$\text{log}_2\vert K(\varepsilon_2)_{_2}\vert$.

In the further considerations we shall always assume that, if $K$
is of the first kind with respect to 2, then $m$ is the constant
of $K$ with respect to 2.

Let $K\not= K(\varepsilon_2)$. Then
$(K(\varepsilon_2):K)=2$. For every element
$\alpha\in K(\varepsilon_2)$ one defines the norm $N(\alpha)$ over $K$
by the equality $N(\alpha)=\alpha\bar\alpha$. The norm is a
multiplicative function. Obviously $N(\alpha)=\alpha^2$ if and only if
$\alpha\in K$. We shall use these elementary properties of the norm
in the next lemma.

\proclaim{Lemma 2} Let $K\not= K(\varepsilon_2)$ and let
$\alpha^2=c\varepsilon_t$, where
$\alpha\in K(\varepsilon_2)^{\ast}$,
$\varepsilon_t\in K(\varepsilon_2)^{\ast}$
and $N(\alpha)=c\in K^{\ast}$. Then

{\rm 1)} If $K$ is of type C, then $c\in K^{\ast^2}$.

{\rm 2)} If $K$ is of type D or E and $0\leq t\leq m-2$, then
$c\in K^{\ast^2}$.

{\rm 3)} If $K$ is of the first kind with respect to $2$ and $t=m-1$, then
$$
c\in
\cases K^{\ast^2},& \text{ if $K$ is of type D,}\\
-K^{\ast^2},& \text{ if $K$ is of type E.}
\endcases
$$

{\rm 4)} If $K$ is of the first kind with respect to $2$ and
$t=m$, then $K$ is of type D and
$c\in(\varepsilon_m+\varepsilon_m^{-1}+2)K^{\ast^2}$.
\endproclaim

\demo{Proof} 1) Let $K$ be of type C. Then $K$ is of the second kind
with respect to 2 and $N(\varepsilon_{t+1})=1$ for every nonnegative
integer $t$. Let $\beta=\alpha\varepsilon_{t+1}^{-1}$. We obtain
immediately $c=\beta^2$. Besides, $N(\beta)=N(\alpha)=c=\beta^2$.
Hence $\beta\in K^{\ast}$ and $c\in K^{\ast^2}$.

2) Let $K$ be of type D or E. Then $N(\varepsilon_{t+1})=1$
for every integer $t\in [0,m-2]$ and the proof goes in the same way as
in the case 1).

3) If $t=m-1$ and $K$ is of type D, then
$N(\varepsilon_{t+1})=N(\varepsilon_m)=1$ and the proof is the same as
in the previous two cases. Now, let $t=m-1$ and let $K$ be of type E.
Then $N(\varepsilon_{t+1})=N(\varepsilon_m)=-1$. This time we set
$\beta=\alpha\varepsilon_2\varepsilon_m^{-1}$. From here we derive that
$c=-\beta^2$. Moreover, $N(\varepsilon_2)=1$ and we have
$N(\beta)=-N(\alpha)=-c=\beta^2$. Hence $\beta\in K^{\ast}$ and
$c\in -K^{\ast^2}$.

4) Let $K$ be of the first kind with respect to 2 and let $t=m$. Then
from the equalities $\alpha^2=c\varepsilon_m$ and $N(\alpha)=c$
given in the condition of the lemma
we obtain $N(\varepsilon_m)=1$. Hence $K$
is of type D. Let $\beta=\alpha(1+\varepsilon_m)^{-1}$. This is possible
because $\varepsilon_m\not=-1$ and immediately implies that
$c=(\varepsilon_m+\varepsilon_m^{-1}+2)\beta^2$. Since $K$ is of type D,
we obtain that $\bar\varepsilon_m=\varepsilon_m^{-1}$. Hence we have
$$
N(\beta)=N(\alpha)(1+\varepsilon_m)^{-1}(1+\bar\varepsilon_m)^{-1}=
$$
$$
=c(1+\varepsilon_m)^{-1}(1+\varepsilon_m^{-1})^{-1}=
c(\varepsilon_m+\varepsilon_m^{-1}+2)^{-1}=\beta^2.
$$
Therefore $\beta\in K^{\ast}$ and
$c\in(\varepsilon_m+\varepsilon_m^{-1}+2)K^{\ast^2}$, which completes
the proof of the lemma.
\enddemo

With every nonnegative integer $s$ we relate
$K_s=K^{\ast}\cap K(\varepsilon_2)^{\ast^{2^s}}$. Clearly $K_s$ is a
multiplicative group which is a subgroup of $K^{\ast}$. One sees immediately
that if $K=K(\varepsilon_2)$, then $K_s=K^{\ast^{2^s}}$. In the next
lemma we give the inner structure of the group $K_s$ when
$K\not= K(\varepsilon_2)$.

\proclaim{Lemma 3} Let $K\not= K(\varepsilon_2)$ and let $s$ be
a nonnegative integer. Then for the group $K_s$ we have:

{\rm 1)} $K_s=K^{\ast}$ for $s=0$.

{\rm 2)} If $K$ is of type C and $s\geq 1$, then
$K_s=\langle -1\rangle\times K^{\ast^{2^s}}$.

{\rm 3)} If $K$ is of type D or E and $1\leq s\leq m-1$, then
$K_s=\langle -1\rangle\times K^{\ast^{2^s}}$.

{\rm 4)} If $K$ is of type D and $s\geq m$, then
$K_s=\langle(1+\varepsilon_m)^{2^s}\rangle K^{\ast^{2^s}}$.

{\rm 5)} If $K$ is of type E and $s\geq m$, then
$K_s=K^{\ast^{2^s}}$.
\endproclaim

\demo{Proof} The case 1) is trivial. In the other cases it is easy
to see that $K_s$ contains the group from the corresponding case.
In the cases 2) and 3) this follows from the formula
$\varepsilon_{s+1}^{2^s}=-1$ and $\varepsilon_{s+1}\in K(\varepsilon_2)$.
For the case 4) the assertion follows from the inclusion
$(\varepsilon_m+1)^{2^s}\in K(\varepsilon_2)^{\ast^{2^s}}$
and in the case 5) it is obvious.

We shall prove the inverse inclusion. For this purpose, let $a\in K_s$.
Then $a\in K^{\ast}$ and
$$
a=\alpha^{2^s},
\leqno(1)
$$
where $\alpha\in K(\varepsilon_2)^{\ast}$. Let us denote
$N(\alpha)=c\in K^{\ast}$. Applying the norm to the equality (1)
we obtain $a^2=c^{2^s}$ and the squaring of (1) gives
$a^2=\alpha^{2^{s+1}}$. Hence
$$
\alpha^{2^{s+1}}=c^{2^s},
\leqno(2)
$$
and taking the root of (2) we obtain
$$
\alpha^2=c\varepsilon_t,
\leqno(3)
$$
where $t\leq s$. If $K$ is of type C, then $t$ is a nonnegative integer
satisfying only the inequality $t\leq s$. If $K$ is of type D or E,
then $t$ is again a nonnegative integer and satisfies
$t\leq\text{min}(s,m)$. Now we can apply Lemma 2 to the equality (3).

Case 2). If $K$ is of type C and $s\geq 1$, then Case 1) of Lemma 2
gives $c=c_1^2$, $c_1\in K^{\ast}$. From here and from (2) we derive
$\alpha^{2^{s+1}}=c_1^{2^{s+1}}$ and $\alpha=c_1\varepsilon_{s+1}^l$
for some integer $l$. Then from (1) we obtain
$a=c_1^{2^s}(-1)^l\in\langle-1\rangle \times K^{\ast^{2^s}}$.

Case 3). If $K$ is of type D or E and $1\leq s\leq m-1$, then
$0\leq t\leq m-1$. Then by Lemma 2, Case 2) or Case 3), we have
$c=c_1^2$ or $c=-c_1^2$, where $c_1\in K^{\ast}$. From here and
from (2) we obtain that $\alpha^{2^{s+1}}=c_1^{2^{s+1}}$ and,
as in the case 2) it follows
$a\in\langle -1\rangle\times K^{\ast^{2^s}}$.

Case 4). If $K$ is of type D and $s\geq m$, then $0\leq t\leq m$. If
$0\leq t\leq m-1$, then Case 2) of Lemma 2 gives $c=c_1^2$,
$c_1\in K^{\ast}$. Together with (2) this implies
$\alpha=c_1\varepsilon_m^l$, $l\in {\Bbb Z}$. Then from (1) we have
$$
a=c_1^{2^s}\varepsilon_m^{l2^s}\in K^{\ast^{2^s}}\subseteq
\langle(1+\varepsilon_m)^{2^s}\rangle K^{\ast^{2^s}}.
$$
Now, if $t=m$, then by Case 4) of Lemma 2 we have
$$
c=(\varepsilon_m+\varepsilon_m^{-1}+2)c_1^2=
\varepsilon_m^{-1}(1+\varepsilon_m)^2c_1^2,\, c_1\in K^{\ast}.
$$
From here and from (2) it follows
$\alpha=c_1(1+\varepsilon_m)\varepsilon_m^l$, $l\in{\Bbb Z}$.
Then (2) gives
$$
a=(1+\varepsilon_m)^{2^s}c_1^{2^s}\in
\langle(1+\varepsilon_m)^{2^s}\rangle K^{\ast^{2^s}}.
$$

Case 5). If $K$ is of type E and $s\geq m$, then $0\leq t\leq m$.
The case $t=m$ is impossible by Case 4) of Lemma 2 and hence
$0\leq t\leq m-1$. Then Lemma 2 Case 2) or Case 3) implies that
$c=c_1^2$ or $c=-c_1^2$, $c_1\in K^{\ast}$.
The minus sign does not influence when replacing $c$ in (2).
Now, proceeding as in Case 4) we obtain that
$a=c_1^{2^s}\in K^{\ast^{2^s}}$. The lemma is established.
\enddemo

\head 3. Minimal idempotents of twisted group algebras \endhead

Let us denote by $s=H_n(a)$ the greatest integer from the interval
$[0,n]$ such that $a\in K_s$.

In the next theorem we use the following elementary fact. Let $A$ be a
commutative $K$-algebra and let $B$ be a subalgebra of $A$. Then, if every
minimal idempotent of $B$ is minimal also in $A$, then the set of
minimal idempotents of $A$ coincides with the set of minimal idempotents
of $B$.

\proclaim {Theorem 1} Let $K=K(\varepsilon_2)$, let
$G=\langle g\rangle$ be a cyclic group of order $2^n$ and
let $K_tG$ be the twisted group algebra of $G$ over $K$
defined by the equality $\bar g^{^{2^n}}=a\in K^{\ast}$.
Let $H_n(a)=s$ and $a=b^{2^s}$, $b\in K^{\ast}$. Let us denote
$h=\bar g^{^{2^{n-s}}}b^{-1}$. Then

{\rm 1)} $H=\langle h\rangle$ is a group of order $2^s$ and the
elements of this group are linearly independent over $K$, i.e.
the group algebra $KH$ is a subalgebra of $K_tG$.

{\rm 2)} The set of minimal idempotents of $K_tG$ coincides with the set
of minimal idempotents of $KH$.
\endproclaim

\demo{Proof} 1) We have $h^{2^s}=\bar g^{^{2^n}}b^{-2^s}=aa^{-1}=1$. If
we assume that $h^{2^{s-1}}=1$, then this implies that
$\bar g^{^{2^{n-1}}}=b^{2^{s-1}}$ which contradicts with the linear
independence of the elements $1,\bar g,\bar g^{^2},\ldots,\bar g^{^{2^n-1}}$.
Hence $h$ is of order $2^s$. If we assume that there exists a linear
dependence
$$
a_0+a_1h+\ldots+a_{2^s-1}h^{2^s-1}=0
$$
and replace $h=\bar g^{^{2^{n-s}}}b^{-1}$, then this immediately brings
to the linear dependence of the powers of $\bar g$. Hence $H$ is
linearly independent over $K$ and $KH$ is a subalgebra of $K_tG$.

2) Let $e$ be a minimal idempotent of $KH$.
Then $KHe\cong K(\varepsilon_t)$, where $0\leq t\leq s$. If
$\varepsilon_t\in K$, then $he=\varepsilon_t^ie$ for some integer $i$.
Then $(\bar ge)^{2^{n-s}}=(be)(he)=(b\varepsilon_t^i)e$ and this means that
$\bar g$ is a root of the polynomial
$f(x)=x^{2^{n-s}}-b\varepsilon_t^i$. This polynomial is irreducible
over $K$ since in the opposite case by Lemma B we have
$b\varepsilon_t^i=b_1^2$, $b_1\in K$. Then
$a=b^{2^s}=(b\varepsilon_t^i)^{2^s}=b_1^{2^{s+1}}$ and for $s\not= n$
we obtain a contradiction with $H_n(a)=s$. (For $s=n$ the case is clear
because then $K_tG=KH$.) Therefore $\bar g$ is a root of
an irreducible over $K$ polynomial. But $K_tG$ is generated by $\bar ge$ as a
$K$-algebra. Hence $K_tG$ is a field and $e$ is a minimal idempotent of
$K_tG$.

If $\varepsilon_t\not\in K$, then Lemma 1 gives that $K$ is of type B
with constant $m<s$ and $m<t\leq s$. Then
$(he)^{2^{t-m}}=\varepsilon_m^ie$ for some odd integer $i$. Hence we
obtain
$$
(\bar ge)^{2^{n-s+t-m}}=
(b^{2^{t-m}}e)(he)^{2^{t-m}}=b^{2^{t-m}}\varepsilon_m^ie.
$$
Therefore $\bar ge$ is a root of the polynomial
$f(x)=x^{2^{n-s+t-m}}-b^{2^{t-m}}\varepsilon_m^i$ which according
to Lemma B is also irreducible over $K$ since $t-m\geq 1$ and
$\varepsilon_m^i\not\in K^{\ast^2}$. This implies that $(K_tG)e$
is a field which proves the theorem.
\enddemo

Using this theorem we shall find the minimal idempotents of $K_tG$
when $K=K(\varepsilon_2)$.

\proclaim {Theorem 2} Let $K=K(\varepsilon_2)$, let
$G=\langle g\rangle$ be a cyclic group of order $2^n$ and let
$K_tG$ be defined by $\bar g^{^{2^n}}=a$. Let $H_n(a)=s$ and
$a=b^{2^s}$, $b\in K^{\ast}$. Then the minimal idempotents
$e_i$ and $e_{ri}$ of $K_tG$ are the following:

{\rm 1)} If at least one of the following conditions holds:

{\rm a)} $K$ is of type A,

{\rm b)} $K$ is of type B and $0\leq s\leq m$ ($m$ is the constant of $K$),

\noindent then
$$
e_i={1\over 2^s}\sum_{j=0}^{2^s-1}
\left(b^{-1}\varepsilon_s^{-i}\bar g^{^{2^{n-s}}}\right)^j,\,
i=0,1,2,\ldots,2^s-1.
\leqno(4)
$$

{\rm 2)} If $K$ is of type B and $m<s\leq n$, then
$$
e_i={1\over 2^s}\sum_{j=0}^{2^s-1}
\left(b^{-1}\varepsilon_m^{-i}\bar g^{^{2^{n-s}}}\right)^j,\,
i=0,1,2,\ldots,2^m-1,
\leqno(5)
$$
and
$$
e_{ri}={1\over 2^{s-r}}\sum_{j=0}^{2^{s-r}-1}
\left(b^{-2^r}\varepsilon_m^{-1}\varepsilon_{m-1}^{-i}
\bar g^{^{2^{n-s+r}}}\right)^j,
\leqno(6)
$$
$r=1,2,\ldots,s-m$, $i=0,1,2,\ldots,2^{m-1}-1$.
\endproclaim

\demo{Proof} Let $H=\langle h\rangle$, where
$h=\bar g^{^{2^{n-s}}}b^{-1}$.
In the case 1) we have $\varepsilon_s\in K$ and then the
minimal idempotents of $KH$
are
$$
e_i={1\over 2^s}\sum_{j=0}^{2^s-1}
\left(\varepsilon_s^{-i}h\right)^j,\,i=0,1,2,\ldots,2^s-1.
$$
In the case 2) the minimal idempotents of $KH$ are
$$
e_i={1\over 2^s}\sum_{j=0}^{2^s-1}
\left(\varepsilon_m^{-i}h\right)^j,\,i=0,1,2,\ldots,2^m-1,
$$
and
$$
e_{ri}={1\over 2^{s-r}}\sum_{j=0}^{2^{s-r}-1}
\left(\varepsilon_m^{-1}\varepsilon_{m-1}^{-i}h^{2^r}\right)^j,\,
r=1,2,\ldots,s-m, i=0,1,2,\ldots,2^{m-1}-1.
$$
Replacing in the above formulas $h=\bar g^{^{2^{n-s}}}b^{-1}$ we obtain
(4), (5) and (6). Then, by Theorem 1, the idempotents (4), (5) and (6)
give a complete list of the minimal idempotents of $K_tG$.
With this the theorem is established.
\enddemo

\proclaim{Theorem 3} Let $K\not= K(\varepsilon_2)$, let
$G=\langle g\rangle$ be a cyclic group of order $2^n$, let $K_tG$ be
defined by $\bar g^{^{2^n}}=a\in K^{\ast}$ and let $H_n(a)=s$. Then
we have for the minimal idempotents of $K_tG$:

{\rm 1)} If $s=0$, then $K_tG$ has a unique minimal idempotent which is
equal to the unity element.

{\rm 2)} If at least one of the following conditions is fulfilled:

{\rm a)} $K$ is of type C, $s\geq 1$ and $a\in K^{\ast^{2^s}}$,

{\rm b)} $K$ is of type D or E, $1\leq s\leq m-1$ and $a\in K^{\ast^{2^s}}$,

\noindent then the minimal idempotents are
$$
e_0={1\over 2^s}\sum_{j=0}^{2^s-1}(b^{-1}\bar g^{^{2^{n-s}}})^j,\,
e_{2^{s-1}}={1\over 2^s}\sum_{j=0}^{2^s-1}(-1)^j(b^{-1}\bar g^{^{2^{n-s}}})^j,
$$
$$
\leqno(7)
$$
$$
e_i={1\over 2^s}\sum_{j=0}^{2^s-1}
(\varepsilon_s^{ij}+\varepsilon_s^{-ij})
(b^{-1}\bar g^{^{2^{n-s}}})^j,\,i=1,2,\ldots,2^{s-1}-1,
$$
where $a=b^{2^s}$, $b\in K^{\ast}$.

{\rm 3)} If $K$ is of type D or E, $s\geq m$ and $a\in K^{\ast^{2^s}}$,
then the minimal idempotents are
$$
e_0={1\over 2^s}\sum_{j=0}^{2^s-1}(b^{-1}\bar g^{^{2^{n-s}}})^j,\,
e_{2^{m-2}}={1\over 2^s}\sum_{j=0}^{2^s-1}(-1)^j(b^{-1}\bar g^{^{2^{n-s}}})^j,
$$
$$
e_i={1\over 2^s}\sum_{j=0}^{2^s-1}
(\varepsilon_{m-1}^{ij}+\varepsilon_{m-1}^{-ij})
(b^{-1}\bar g^{^{2^{n-s}}})^j,\,i=1,2,\ldots,2^{m-2}-1,
\leqno(8)
$$
$$
e_{ri}={1\over 2^{s-r}}\sum_{j=0}^{2^{s-r}-1}
(\varepsilon_m^{-j}\varepsilon_{m-2}^{-ij}+
\lambda^{-j}\varepsilon_m^j\varepsilon_{m-2}^{ij})
(b^{-2^r}\bar g^{^{2^{n-s+r}}})^j,
$$
$$
r=0,1,2,\ldots,s-m,\,i=0,1,2,\ldots,2^{m-2}-1,
$$
where $a=b^{2^s}$, $b\in K^{\ast}$ and
$\lambda=\cases 1,& \text{ if $K$ is of type D,}\\
-1,& \text{ if $K$ is of type E.}
\endcases$.
\bigskip

{\rm 4)} If at least one of the following conditions is fulfilled:

{\rm a)} $K$ is of type C, $s\geq 1$ and $a\in -K^{\ast^{2^s}}$,

{\rm b)} $K$ is of type D or E, $1\leq s\leq m-1$ and $a\in -K^{\ast^{2^s}}$,

\noindent then the minimal idempotents are
$$
e_i={1\over 2^s}\sum_{j=0}^{2^s-1}
(\varepsilon_{s+1}^{-j}\varepsilon_{s-1}^{-ij}+
\lambda^{-j}\varepsilon_{s+1}^j\varepsilon_{s-1}^{ij}
(b^{-1}\bar g^{^{2^{n-s}}})^j,\,i=1,2,\ldots,2^{s-1}-1,
\leqno(9)
$$
$a=-b^{2^s}$, $b\in K^{\ast}$ and
$$
\lambda=\cases
-1,&\text{ if $K$ is of type E and $s=m-1$,}\\
1,& \text{ in all other cases.}
\endcases
$$

{\rm 5)} If $K$ is of type D, $s\geq m$ and
$a\in(1+\varepsilon_m)^{2^s}K^{\ast^{2^s}}$, then
the minimal idempotents are
$$
e_i={1\over 2^s}\sum_{j=0}^{2^s-1}(1+\varepsilon_m)^{-j}
(\varepsilon_{m-1}^{-ij}+
\varepsilon_m^j\varepsilon_{m-1}^{ij})
(b^{-1}\bar g^{^{2^{n-s}}})^j,\,i=0,1,2,\ldots,2^{m-1}-1,
$$
$$
e_{1,2^{m-1}-1}={1\over 2^{s-1}}\sum_{j=0}^{2^{s-1}-1}
(2+\varepsilon_m+\varepsilon_m^{-1})^{-j}(b^{-2}\bar g^{^{2^{n-s+1}}})^j,\,
s\geq m+1,
$$
$$
e_{1,2^{m-2}-1}={1\over 2^{s-1}}\sum_{j=0}^{2^{s-1}-1}(-1)^j
(2+\varepsilon_m+\varepsilon_m^{-1})^{-j}(b^{-2}\bar g^{^{2^{n-s+1}}})^j,\,
s\geq m+1,
\leqno(10)
$$
$$
e_{1,i}={1\over 2^{s-1}}\sum_{j=0}^{2^{s-1}-1}
(2+\varepsilon_m+\varepsilon_m^{-1})^{-j}
(\varepsilon_{m-1}^{-j(1+i)}+
\varepsilon_{m-1}^{j(1+i)})
(b^{-2}\bar g^{^{2^{n-s+1}}})^j,
$$
$$
s\geq m+1,\,m\geq 3,\,i=0,1,2,\ldots,2^{m-2}-2,
$$
$$
e_{ri}={1\over 2^{s-r}}\sum_{j=0}^{2^{s-r}-1}
(1+\varepsilon_m)^{-j2^r}
(\varepsilon_m^{-j}\varepsilon_{m-2}^{-ij}+
\varepsilon_m^j\varepsilon_{m-2}^{j(i+2^{r-2})})
(b^{-2^r}\bar g^{^{2^{n-s+r}}})^j,
$$
$$
r=2,3,\ldots,s-m,\,i=0,1,2,\ldots,2^{m-2}-1,\,s\geq m+2
$$
and for all these formulas we have
$a=b^{2^s}(1+\varepsilon_m)^{2^s}$, $b\in K^{\ast}$.
\endproclaim

\demo{Proof} The case $s=0$ is clear because, by Theorem A,
$K_tG$ is a field. In the other cases we first find the minimal
idempotents of $K(\varepsilon_2)_tG$ according to Theorem 2
replacing $b$ with $y\in K(\varepsilon_2)$, where $y$ is an arbitrary root
of the equation $y^{2^s}=a$. We determine $y$ by Lemma 3, where we can
take a suitably chosen root of this equation. Then we group the
idempotents of $K(\varepsilon_2)_tG$ in pairs of conjugates and add
the elements from the pairs. In this way we obtain the minimal
idempotents of $K_tG$. One should have in mind that some of the
idempotents may be conjugate to themselves. Then they will be also
idempotents of $K_tG$. Every idempotent of $K(\varepsilon_2)_tG$
is indexed by a single index $i$ or by a double index $ri$. Then we
index the conjugate respectively with $i'$ and $ri'$. There exists a
relation of the kind $i'+i=\text{const}$ between the indices $i'$ and $i$.
An idempotent $e$ is self-conjugate only if $i'=i$. This is possible
only if this constant is even. When the constant is odd, in every class
of conjugate idempotents there is one even and one odd index. Then
the even indices can serve as a complete system of representatives of
the classes of conjugate idempotents.

For the case 2) we use formula (4) from Theorem 2, where
$y=b\in K^{\ast}$. The formula for the conjugacy is $i'+i=2^s$ if
$i\not=0$ and $i'=0$ for $i=0$, since
$\overline{\varepsilon_s^i}=\varepsilon_s^{-i}$. For $i=0$ and $i=2^{s-1}$
the idempotents are self-conjugate. For $i=1,2,\ldots,2^{s-1}-1$
we have a complete system of representatives of the classes of conjugate
idempotents. In this way we obtain the formula (7).

In the case 3) we use the formula (4) for $s=m$ and the formulas
(5) and (6) from theorem 2 when $y=b\in K^{\ast}$. For $s=m$ the formula
(4) is identical with (5) and we shall use only (5) and (6),
allowing in (5) also $s=m$. Moreover, we present the formulas (5) and (6)
in another way, assuming that the index $i$ in (5) has only even values
and the odd values are transferred in (6) with $r=0$. Then we introduce
a new index $i=0,1,2,\ldots,2^{m-1}-1$ in (5), and change $m$ with $m-1$.
In (6) we allow also $r=0$.Using the transformed versions of the
formulas (5) and (6) we look for the conjugate idempotents. The
formula for conjugating of the idempotents from (5) is
$i'+i=2^{m-1}$ for $i\not=0$ and $i'=0$ for $i=0$. Then we have two
self-conjugate idempotents for $i=0$ and $i=2^{m-2}$. For each of
the other idempotents for $i=1,2,\ldots,2^{m-2}-1$, we have one class
conjugate idempotents. In this way we obtain the first three formulas
(without the last) from (8).

The formula for conjugating of the idempotents from (6) when $K$
is of type D is $i'+i=2^{m-1}-1$, and when $K$ is of type E, it is
$$
i'+i=\cases
2^{m-2}-1,& \text{ if } 0\leq i\leq 2^{m-2}-1,\\
3.2^{m-2}-1,& \text{ if } 2^{m-2}\leq i\leq 2^{m-1}-1,
\endcases
\leqno(11)
$$
because in the first case we have
$\bar\varepsilon_m=\varepsilon_m^{-1}$ and in the second case
$\bar\varepsilon_m=-\varepsilon_m^{-1}$. We see from these formulas that
$i'\not=i$ and hence there are no self-conjugate idempotents.
When $i$ runs on the set of all even integers from $0$ to $2^{m-1}-1$,
we obtain a complete system of representatives of the classes of
conjugate idempotents. Then, introducing a new index $i$ equal to the half
of the previous, changing $m-1$ with $m-2$ and taking into account that
$\bar\varepsilon_m=\lambda\varepsilon_m^{-1}$, we obtain also the
last formula from (8).

In the case 4) we use the formula (4) of Theorem 2 where we have
$y=b\varepsilon_{s+1}$ instead of $b$. When $K$ is of type E and
$s=m-1$ the conjugacy formula is (11) and in all other cases it is
$i'+i=2^s-1$. Hence $i'\not=i$ and therefore we have no self-conjugate
idempotents. In order to have a complete system of representatives
of the classes of conjugate idempotents we can use the even integers from
$0$ to $2^s-1$ and, as in the case 3), changing the index $i$ we obtain
the formula (9).

In the case 5) we use the formulas (5) and (6) from Theorem 2 where
instead of $b$ we have $y=b(1+\varepsilon_m)$ and in the formula (5)
we allow also $s=m$. The conjugacy formula for the idempotents from (5)
is $i'+i=2^m-1$ and, as before, we see that there are no self-conjugate
idempotents and the even integers from $0$ to $2^m-1$ serve for a complete
system of representatives of the classes of conjugate idempotents.
Again, changing the index $i$ we obtain the first formula from (10).
For the idempotents from (6) and for $r=1$ the conjugacy formula is
$i'+i=2^{m-1}-2$ when $i\not=2^{m-1}-1$ and $i'=2^{m-1}-1$ for
$i=2^{m-1}-1$. From here we obtain that for $i=2^{m-1}-1$ and
$i=2^{m-2}-1$ the idempotents are self-conjugate. These are the
idempotents $e_{1,2^{m-1}-1}$ and $e_{1,2^{m-2}-1}$. They exist for
$s\geq m+1$ only because this inequality is required in (6). The other
idempotents are grouped in pairs of conjugates when $i$ has values
from $0$ to $2^{m-2}-2$. The condition for the existence of these
idempotents is not only $s\geq m+1$. We need also $m\geq 3$ because
for $m=2$ the only idempotents are the self-conjugate ones. In this
way we obtain all other formulas from (10) except the last one.
For the idempotents from (6) for $r\geq 2$ the conjugacy formula is
$$
i'+i=\cases
2^{m-1}-2^{r-1}-1,& \text{ if $r<m$ and $0\leq i\leq 2^{m-1}-2^{r-1}-1$}\\
2^m-2^{r-1}-1,& \text{ if $r<m$ and
$2^{m-1}-2^{r-1}\leq i\leq 2^{m-1}-1$}\\
2^{m-1}-1,& \text{ if $r\geq m$.}
\endcases
$$
We see from this formula that the sum $i'+i$ is odd and hence there are
no self-conjugate idempotents. Then we can choose for a complete system
of representatives of the classes of conjugate idempotents using
the even integers from the interval $[0,2^{m-1}-1]$ and after a change
of the index $i$ we obtain the last formula from (10). The existence
condition for these idempotents is $s\geq m+2$ since $s-m\geq r\geq 2$.
The theorem is established.
\enddemo

\enddocument
\end